\documentclass[a4paper,10pt]{article}

\usepackage{fullpage}

\usepackage{natbib}

\usepackage{lmodern}
\usepackage{hyperref}
\usepackage{amscd}
\usepackage{mathtools}
\usepackage{amssymb}
\mathtoolsset{showonlyrefs=true}
\usepackage{dsfont}

\def\1{\mathds{1}}
\def\P{\mathbb{P}} 
\def\E{\mathbb{E}} 
\def\Re{{\mathds{R}}}
\def\G{{\cal G}} 
\def\T{{\cal T}}
\def\R{{\cal R}}
\def\sw{S}
\def\switch{\mbox{switch}}
\def\Vmax{V_{\max}}

\def\hpi{{Howard's PI}}
\def\spi{{Simplex-PI}}

\def\titre{Improved and Generalized Upper Bounds on the Complexity of Policy Iteration}

\usepackage{amsthm}
\usepackage{amscd}
\usepackage{amsfonts}

\newcommand{\appdx}[1]{Section~\ref{#1}}

\def\beginproof{\begin{proof}}
\def\endproofb{\end{proof}}

\def\fornips{For clarity, all proofs are deferred to the later sections.}

\newtheorem{theorem}{Theorem}
\newtheorem{lemma}{Lemma}
\newtheorem{definition}{Definition}
\newtheorem{corollary}{Corollary}

\newtheorem{assumption}{Assumption}
\newtheorem{remark}{Remark}

\title{\titre}

\author{Bruno Scherrer \\
INRIA Nancy Grand Est, Team MAIA\\
bruno.scherrer@inria.fr}

\begin{document}

\maketitle

\begin{abstract}%
Given a Markov Decision Process (MDP) with $n$ states and a total number $m$ of actions, we study the number of iterations needed by Policy Iteration (PI) algorithms to converge to the optimal $\gamma$-discounted policy. We consider two variations of PI: \hpi{ }that changes the actions in all states with a positive advantage, and \spi{ }that only changes the action in the state with maximal advantage.
We show that \hpi{ }terminates after at most 
$
 O \left( \frac{m}{1-\gamma} \log \left( \frac{1}{1-\gamma} \right)\right)
$
iterations, improving by a factor $O(\log n)$ a result by Hansen et al (2013), while \spi{ }terminates after at most
$
O \left(  \frac{n m}{1-\gamma} \log \left( \frac{1}{1-\gamma} \right)\right)
$
iterations, improving by a factor $O(\log n)$ a result by Ye (2011).
Under some structural properties of the MDP, we then consider bounds that are independent of the discount factor~$\gamma$: quantities of interest are bounds $\tau_t$ and $\tau_r$---uniform on all states and policies---respectively on the \emph{expected time spent in transient states} and \emph{the inverse of the frequency of visits in recurrent states} given that the process starts from the uniform distribution.
Indeed, we show that \spi{ }terminates after at most
$
\tilde O \left( n^3 m^2 \tau_t \tau_r \right)
$
iterations. This extends a recent result for deterministic MDPs by Post \& Ye (2013), in which $\tau_t \le 1$ and $\tau_r \le n$; in particular it shows that \spi{ }is strongly polynomial for a much larger class of MDPs. We explain why similar results seem hard to derive for \hpi. Finally, under the additional (restrictive) assumption that the state space is partitioned in two sets, respectively states that are transient and recurrent for all policies, we show that both \hpi{ }and \spi{ }terminate after at most 
$
\tilde O(m (n^2 \tau_t+ n \tau_r))$ iterations.

\end{abstract}


\section{Introduction}

We consider a discrete-time dynamic system whose state transition
depends on a control, where the {\bf state space} $X$ is of finite size $n$. When at state $i \in \{1,..,n\}$, the action is chosen from a set of admissible actions $A_i \subset A$, where the {\bf action space}
$A$ is of finite size $m$, such that $(A_i)_{1 \le i \le n}$ form a partition of $A$. The action $a \in A_i$ specifies the {\bf transition probability} $p_{ij}(a)=\P(i_{t+1}=j | i_t=i, a_t=a)$
to the next state $j$. At each transition, the system is given a reward $r(i,a,j) \in \mathbb{R}$ where $r$ is the instantaneous {\bf reward function}. In this context, we look for a stationary deterministic policy\footnote{Restricting our attention to stationary deterministic policies is not a limitation. Indeed, for the optimality criterion to be defined soon, it can be shown that there exists at least one
stationary deterministic policy that is optimal
\citep{puterman}.}, that is a
function $\pi:X \rightarrow A$ that maps states into
admissible actions (for all $i$, $\pi(i) \in A_i$) that maximizes the expected discounted sum of rewards from any state~$i$, called  the {\bf value
of policy $\pi$} at state~$i$:
\begin{equation}
\label{bellvdef}
v_\pi(i):= \E \left[ \left. \sum_{k=0}^{\infty}\gamma^k r(i_k,a_k,i_{k+1}) \right| i_0=i, ~\forall k \ge 0,~a_k=\pi(i_k),~i_{k+1}\sim \P(\cdot | i_k, a_k) \right],
\end{equation}
where $\gamma \in (0,1)$ is a discount factor. The tuple $\langle X, (A_i)_{i \in X}, p, r, \gamma \rangle$ is called a {\bf Markov Decision Process (MDP)} \citep{puterman,ndp}, and the associated problem is known as {\bf stochastic optimal control.}

The {\bf optimal value} starting from state $i$ is defined as
$$
v_*(i):= \max_\pi v_\pi(i).
$$
For any policy $\pi$, we write $P_\pi$ for the $n \times n$ stochastic matrix whose elements are
$p_{ij}(\pi(i))$, and $r_\pi$ for the vector whose components are
$\sum_j p_{ij}(\pi(i))r(i,\pi(i),j)$. The value functions $v_\pi$ and $v_*$ can
be seen as vectors on $X$.  It is well known that $v_\pi$ is the solution of the
following Bellman equation:
$$
v_\pi = r_\pi + \gamma P_\pi v_\pi,
$$
that is $v_\pi$ is a fixed point of the affine
operator $T_\pi:v \mapsto r_\pi + \gamma P_\pi v$. 
It is also well known that $v_*$ satisfies the
following Bellman equation:
$$
v_* = \max_\pi (r_\pi + \gamma P_\pi v_*) = \max_\pi T_\pi v_*
$$
where the \mbox{max} operator is taken componentwise.
In other words, $v_*$ is a fixed point of the nonlinear operator $T:v \mapsto \max_{\pi}T_\pi v$. For any value vector
$v$, we say that a policy $\pi$ is {\bf greedy with respect to the value $v$} if it satisfies:
$$
\pi \in \arg\max_{\pi'} T_{\pi'} v
$$
or equivalently $T_\pi v = T v$. With some slight abuse of notation, we write
$\G(v)$ for any policy that is greedy with respect to
$v$. The notions of optimal value function and greedy policies are
fundamental to optimal control because of the following property: any
policy $\pi_*$ that is greedy with respect to the optimal value $v_*$ is an
{\bf optimal policy} and its value $v_{\pi_*}$ is equal to $v_*$.

Let $\pi$ be some policy. For any policy $\pi'$, we consider the quantity 
$$
a_\pi^{\pi'}=T_{\pi'} v_\pi-v_\pi
$$
that measures the difference in value resulting from switching the first action to $\pi'$ with respect to always using $\pi$; we shall call it the {\bf advantage of $\pi'$ with respect to $\pi$}.
Furthermore, we call {\bf maximal advantage with respect to $\pi$} the componentwise best such advantage:
$$
a_\pi=\max_{\pi'} a_\pi^{\pi'} = T v_\pi - v_\pi,
$$
where the second equality follows from the very definition of the Bellman operator $T$. While the advantage $a_\pi^{\pi'}$ may have negative values, the maximal advantage $a_\pi$ has only non-negative values.
We call the {\bf set of switchable states of $\pi$} the set of states for which the maximal advantage with respect to $\pi$ is positive:
$$
\sw_\pi=\{ i,~ a_{\pi}(i)>0 \}.
$$
Assume now that $\pi$ is non-optimal (this implies that $\sw_\pi$ is a non-empty set). For any non-empty subset $Y$ of $\sw_\pi$, we denote $\switch(\pi,Y)$ a policy satisfying: 
$$
\forall i,~ \switch(\pi,Y)(i)=
\left\{
\begin{array}{ll}
\G(v_\pi)(i)  & \mbox{if }i \in Y\\
\pi(i) & \mbox{if }i \not\in Y.
\end{array}
\right. 
$$
The following result is well known (see for instance \citet{puterman}).
\begin{lemma}
\label{piimp}
Let $\pi$ be some non-optimal policy. If $\pi'=\switch(\pi,Y)$ for some non-empty subset  $Y$ of $\sw_\pi$,  then  $v_{\pi'} \ge v_\pi$ and there exists at least one state $i$ such that $v_{\pi'}(i)>v_\pi(i)$.
\end{lemma}
This lemma is the foundation of the well-known iterative procedure, called Policy Iteration (PI), that generates a sequence of policies $(\pi_k)$  as follows.
\begin{align}
\label{policyiteration}
\pi_{k+1} \leftarrow \switch(\pi_k,Y_k)\mbox{ for some set }Y_k\mbox{ such that }\emptyset \subsetneq Y_k \subseteq \sw_{\pi_k}.
\end{align}
The choice for the subsets $Y_k$ leads to different variations of PI. In this paper we will focus on two of them:
\begin{itemize}
\item When for all iterations $k$, $Y_k=\sw_{\pi_k}$, that is one switches the actions in all states with positive advantage with respect to $\pi_k$, the above algorithm is known as \hpi; it can be seen then that $\pi_{k+1} \in \G(v_{\pi_k})$.
\item When for all iterations $k$, $Y_k$ is a singleton containing a state $i_k \in \arg\max_i a_{\pi_k}(i)$, that is if we only switch one action in the state with maximal advantage with respect to $\pi_k$, we will call it \spi\footnote{In this case, PI is equivalent to running the simplex algorithm with the highest-pivot rule on a linear program version of the MDP problem~\citep{ye}.}.   
\end{itemize}
Since it generates a sequence of policies with increasing values, any
variation of PI converges to an optimal policy in a number of
iterations that is smaller than the total number of policies.  In
practice, PI converges in very few iterations. On random MDP
instances, convergence often occurs in time sub-linear in $n$.  The aim
of this paper is to discuss existing and provide new upper bounds on
the number of iterations required by \hpi{ }and \spi{ }that are much
sharper than $m^n$.

In the next sections, we describe some known results---see also \citet{ye} for a recent and comprehensive review---about the number of iterations required by \hpi{ }and \spi, along with some of our original improvements and extensions. \fornips

\section{Bounds with respect to a fixed discount factor $\gamma<1$}

A key observation for both algorithms, that will be central to the results we are about to discuss, is that the sequences they generate satisfy some contraction property\footnote{A sequence of non-negative numbers $(x_k)_{k \ge 0}$ is contracting with coefficient $\alpha$ if and only if for all $k \ge 0$, $x_{k+1} \le \alpha x_k$.}. For any vector $u \in \Re^n$, let $\|u\|_\infty=max_{1 \le i \le n}|u(i)|$ be the max-norm of $u$. Let $\1$ be the vector of which all components are equal to 1. 
\begin{lemma}[e.g. \cite{puterman}, proof in Section~\ref{proofhpicontraction}]
\label{hpicontraction}
The sequence $(\|v_*-v_{\pi_k}\|_\infty)_{k \ge 0}$ built by \hpi{ }is contracting with coefficient $\gamma$. 
\end{lemma}
\begin{lemma}[\citep{ye}, proof in Section~\ref{proofspicontraction}]
\label{spicontraction}
The sequence $(\1^T(v_*-v_{\pi_k}))_{k \ge 0}$ built by \spi{ }is contracting with coefficient $1-\frac{1-\gamma}{n}$.
\end{lemma}
Contraction is a widely known property for \hpi, and it was to our knowledge first proved by \citep{ye} for \spi; we provide simple proofs in this paper for the sake of completeness. 
While the first contraction property is based on the $\|\cdot\|_\infty$-norm, the second can be equivalently expressed in terms of the $\|\cdot\|_1$-norm defined by $\|u\|_1=\sum_{i=1}^n |u(i)|$, since the vectors $v_*-v_{\pi_k}$ are non-negative and thus satisfy $\1^T(v_*-v_{\pi_k})=\|v_*-v_{\pi_k}\|_1$.
Contraction has the following immediate consequence\footnote{For \hpi, we have: $\|v_*-v_{\pi_k}\|_\infty \le \gamma^k \|v_*-v_{\pi_0}\|_\infty \le \gamma^k \Vmax$. Thus, a sufficient condition for $\|v_*-v_{\pi_k}\|_\infty < \epsilon$ is $\gamma^k \Vmax < \epsilon$, which is implied by $k \ge \frac{\log{\frac{\Vmax}{\epsilon}}}{1-\gamma} > \frac{\log{\frac{\Vmax}{\epsilon}}}{\log{\frac 1 \gamma}}$. For \spi, we have  $\|v_*-v_{\pi_k}\|_\infty \le \|v_*-v_{\pi_k}\|_1 \le \left(1-\frac{1-\gamma}{n}\right)^k \|v_*-v_{\pi_0}\|_1 \le \left(1-\frac{1-\gamma}{n}\right)^k n \Vmax$, and the conclusion is similar to that for \hpi.}.
\begin{corollary}
\label{simplecorol}
Let $\Vmax=\frac{\max_\pi\|r_\pi\|_\infty}{1-\gamma}$ be an upper bound on $\|v_\pi\|_\infty$ for all policies $\pi$.
In order to get an $\epsilon$-optimal policy, that is a policy $\pi_k$ satisfying $\|v_*-v_{\pi_k}\|_\infty \le \epsilon$, \hpi{ }requires at most 
$
\left\lceil \frac{\log{\frac{\Vmax}{\epsilon}}}{1-\gamma} \right\rceil
$ iterations,
while \spi{ }requires at most 
$
\left\lceil \frac{n \log{\frac{n \Vmax}{\epsilon}}}{1-\gamma} \right\rceil
$ 
iterations.
\end{corollary} 
These bounds depend on the precision term $\epsilon$, which means that \hpi{ }and \spi{ }are \emph{weakly polynomial} for a fixed discount factor $\gamma$.
 An important breakthrough was recently achieved by \citet{ye} who proved that one can remove the dependency with respect to $\epsilon$, and thus show that \hpi{ }and \spi{ }are \emph{strongly polynomial} for a fixed discount factor $\gamma$. 
\begin{theorem}[\citet{ye}]
\label{thm:ye1}
\spi{ }and \hpi{ }both terminate after at most
$$
(m-n) \left\lceil \frac{n}{1-\gamma} \log \left( \frac{n^2}{1-\gamma} \right)\right\rceil ~=~ O\left(\frac{mn}{1-\gamma}\log\frac{n}{1-\gamma}\right)
$$
iterations.
\end{theorem}
The proof is based on the fact that PI corresponds to the simplex algorithm in a linear programming formulation of the MDP problem. Using a more direct proof---not based on linear programming arguments---\citet{hansen} recently improved the result by a factor $O(n)$ for \hpi.
\begin{theorem}[\citet{hansen}]
\hpi{ }terminates after at most
$$
(m+1) \left\lceil \frac{1}{1-\gamma} \log \left( \frac{n}{1-\gamma} \right)\right\rceil ~=~ O\left( \frac{m}{1-\gamma}\log\frac{n}{1-\gamma}\right)
$$
iterations.
\end{theorem}
Our first results, that are consequences of the contraction property of \hpi{ }(Lemma~\ref{hpicontraction}) are stated in the following theorems.
\begin{theorem}[Proof in Section~\ref{proofdischpi}]
\label{dischpi}
\hpi{ }terminates after at most
$$
(m-n) \left \lceil \frac{1 }{1-\gamma} \log \left( \frac{1}{1-\gamma} \right) \right \rceil ~=~ O\left( \frac{m}{1-\gamma}\log\frac{1}{1-\gamma}\right)
$$
iterations.
\end{theorem}
\begin{theorem}[Proof in Section~\ref{proofdiscspi2}]
\label{discspi2}
\spi{ }terminates after at most
$$
n (m-n) \left( 1 + \frac{2}{1-\gamma}\log \frac{1}{1-\gamma} \right) ~=~ O\left( \frac{mn}{1-\gamma}\log\frac{1}{1-\gamma}\right)
$$
iterations.
\end{theorem}
Both results are a factor $O(\log n)$ better than the previously known
results provided by \citet{hansen} and \citet{ye}. 
These improvements boil down to the use of the $\|\cdot\|_\infty$-norm instead of the $\|\cdot\|_1$-norm at various points of the previous analyses. For \hpi, the resulting arguments constitute a rather simple extension---the overall line of analysis ends up being very simple, and we consequently believe that it could be part of an elementary course on Policy Iteration; note that a similar improvement and analysis was discovered independently by \citet{akian} in a slightly more general setting. For \spi, however, the line of analysis is slightly trickier: it amounts to bound the improvement in value at individual states and requires a bit of bookkeeping; the technique we use is to our knowledge original.

The bound for \spi{ }is a factor $O(n)$
larger than that for \hpi\footnote{Note that it was also the case in Corollary~\ref{simplecorol}.}. However, since one changes only one action per iteration, each
iteration has a complexity that is in a worst-case sense lower by a factor $n$:
the update of the value can be done in time $O(n^2)$ through the
Sherman-Morrisson formula, though in general each iteration of \hpi,
which amounts to compute the value of some policy that may be
arbitrarily different from the previous policy, may require $O(n^3)$
time. Thus, it is remarkable that both algorithms seem to have a
similar complexity.

The linear dependency of the bound for \hpi{ }with respect to $m$ is
optimal~\cite[Chapter~6.4]{hansenphd}. The linear dependency with
respect to $n$ \emph{or} $m$ (separately) is easy to prove for \spi; we conjecture that
\spi's complexity is proportional to $n m$, and thus that our bound is
tight for a fixed discount factor.  The dependency with respect to
the term $\frac{1}{1-\gamma}$ may be improved, but removing it is
impossible for \hpi{ }and very unlikely for \spi. \citet{fearnley}
describes an MDP for which \hpi{ }requires an exponential (in $n$)
number of iterations for $\gamma=1$ and \citet{hollanders} argued that
this holds also when $\gamma$ is in the vicinity of $1$.  Though a
similar result does not seem to exist for \spi{ }in the literature,
\citet{condon} consider four variations of PI that all switch one
action per iteration, and show through specifically designed MDPs that
they may require an exponential (in $n$) number of iterations when
$\gamma=1$.

\section{Bounds for \spi{ }that are independent of $\gamma$}

In this section, we will describe some bounds that do not depend on $\gamma$
but that will be based on some structural properties of the MDP.
On this topic, \citet{ye2} recently showed the following result for deterministic
MDPs.
\begin{theorem}[\citet{ye2}]
\label{detspiye}
If the MDP is deterministic, then \spi{ }terminates after at most 
$O(n^3 m^2 \log^2 n)$
 iterations.
\end{theorem}
Given a policy $\pi$ of a deterministic MDP, states are either on cycles or on paths induced by $\pi$. The core of the proof relies on the following lemmas that altogether show that cycles are created regularly and that significant progress is made every time a new cycle appears; in other words, significant progress is made regularly.
\begin{lemma}[{\citet[Lemma~3.4]{ye2}}]
\label{spidetpart1}
If the MDP is deterministic, after $O(n^2 m \log n)$ iterations, either \spi{ }finishes or a new cycle appears.
\end{lemma}
\begin{lemma}[{\citet[Lemma~3.5]{ye2}}]
\label{spidetpart2}
If the MDP is deterministic, when \spi{ }moves from $\pi$ to $\pi'$ where $\pi'$ involves a new cycle, we have
\begin{align}
\1^T(v_{\pi_*}-v_{\pi'})  \le \left( 1-\frac 1 n \right) \1^T(v_{\pi_*}-v_{\pi}).
\end{align}
\end{lemma}
Indeed, these observations suffice to prove\footnote{This can be done by using arguments similar those for Theorem~\ref{thm:ye1} (see \citet{ye} for details).} that \spi{ }terminates after $O(n^2 m^2 \log \frac{n}{1-\gamma})$. Completely removing the dependency with respect to the discount factor $\gamma$---the term in $O(\log \frac{1}{1-\gamma})$---requires a careful extra work described in \citet{ye2}, which incurs an extra term of order $O(n \log(n))$.

The main result of this section is to show how these results can be extended to a more general setting. While \citet{ye} reason on states that belong to \emph{paths} and \emph{cycles} induced by policies on deterministic MDPs, we shall consider their natural generalization for stochastic MDPs: \emph{transient states} and \emph{recurrent classes} induced by policies. Precisely, we are going to consider bounds---uniform on all policies and states---of the average time 1) spent in transient states and 2) needed to revisit states in recurrent classes.
For any policy $\pi$ and state $i$, denote $\tau^\pi(i,t)$ the expected cumulative time spent in state $i$ until time $t-1$ given than the process starts from the uniform distribution $U$ on~$X$ and takes actions according to $\pi$:
\begin{equation}
\label{def:exp_time}
\tau^\pi(i,t) = \E \left[ \sum_{k=0}^{t-1} {\textbf{1}}_{i_t=i}~|~i_0 \sim U,~ a_t=\pi(i_t)\right] = \sum_{k=0}^{t-1} \P(i_t=i~|~i_0 \sim U,~ a_t=\pi(i_t)),
\end{equation}
where ${\textbf{1}}$ denotes the indicator function.
In addition, consider the vector $\mu^\pi$ on $X$ providing the asymptotic frequency in all states given that policy $\pi$ is used and that the process starts from the uniform distribution $U$:
\begin{equation}
\forall i,~\mu^{\pi}(i) = \lim_{t \rightarrow \infty} \frac{1}{t} \tau^\pi(i,t).
\end{equation}
When the Markov chain induced by $\pi$ is ergodic, and thus admits a unique stationary distribution, $\mu^{\pi}$ is equal to this very stationary distribution.
However, our definition is more general in that policies may induce Markov chains with aperiodicity and/or multiple recurrent classes. 
For any state $i$ that is transient for the Markov chain induced by $\pi$, it is well known that $\lim_{t \rightarrow \infty} \tau^\pi(i,t)<\infty$ and $\mu^\pi(i)=0$. However, for any recurrent state $i$, we know that $\lim_{t \rightarrow \infty} \tau^\pi(i,t)=\infty$ and  $\mu^\pi(i)>0$;  in particular, if $i$ belongs to some recurrent class $\R$, which is reached with probability $q$ from the uniform distribution $U$, then $\frac{q}{\mu^\pi(i)}$ is the expected time between two visits of the state $i$. 

We are now ready to express the structural properties with which we can provide an extension of the analysis of \citet{ye2}.
\begin{definition}
\label{ass1}
Let $\tau_t$ and $\tau_r$ be the smallest finite constants such that for all policies $\pi$ and states $i$,
\begin{align}
\mbox{if $i$ is transient for $\pi$, then~~~} & \lim_{t \rightarrow \infty} \tau^\pi(i,t) \le \tau_t\\
\mbox{else if $i$ is recurrent for $\pi$, then~~~} & \frac{1}{\mu^\pi(i)} \le \tau_r.
\end{align}
\end{definition}
Note that for any finite MDP, these finite constants always exist. With Definition~\ref{ass1} in hand, we can generalize Lemmas~\ref{spidetpart1}-\ref{spidetpart2} as follows.
\begin{lemma}
\label{spistocpart1}
After at most 
$
(m-n) \lceil n^2 \tau_t \log (n^2 \tau_t)\rceil + n \lceil n^2 \tau_t \log (n^2 ) \rceil
$ 
iterations either \spi{ }finishes or a new recurrent class appears.
\end{lemma}
\begin{lemma}
\label{spistocpart2}
When \spi{ }moves from $\pi$ to $\pi'$ where $\pi'$ involves a new recurrent class, we have
\begin{align}
\1^T(v_{\pi_*}-v_{\pi'})  \le \left( 1-\frac{1}{n \tau_r} \right) \1^T(v_{\pi_*}-v_{\pi}).
\end{align}
\end{lemma}
From these generalized observations, we can deduce the following original result.
\begin{theorem}[Proof in \appdx{proofstocspi}]
\label{stocspi}
\spi{ }terminates after at most
$$
\left[ m \lceil n \tau_r \log (n^2 \tau_r) \rceil + (m-n) \lceil n \tau_r \log (n^2 \tau_t) \rceil \right] \left[ (m-n) \lceil n^2 \tau_t \log (n^2 \tau_t)\rceil + n \lceil n^2 \tau_t \log (n^2 ) \rceil \right]
=
\tilde O \left( n^3 m^2 \tau_t \tau_r \right)
$$
iterations.
\end{theorem}
\begin{remark}
This new result extends the result obtained for deterministic MDPs by \citet{ye2} recalled in Theorem~\ref{detspiye}. In the deterministic case, it is easy to see that $\tau_t = 1$ and $\tau_r \le n$. Then, while Lemma~\ref{spistocpart1} is a strict generalization of Lemma~\ref{spidetpart1}, Lemma~\ref{spistocpart2} provides a contraction factor that is slightly weaker than that of Lemma~\ref{spidetpart2}---$\left(1-\frac{1}{n^2}\right)$ instead of  $\left(1-\frac{1}{n}\right)$---, which makes the resulting bound provided in Theorem~\ref{stocspi}  a factor $O(n)$ worse than that of Theorem~\ref{detspiye}. This extra term in the bound is the price paid for making the constant $\tau_r$ (and the vector $\mu_\pi$) independent of the discount factor $\gamma$, that is by presenting our result in a way that only depends on the dynamics of the underlying MDP. An analysis that would strictly generalizes that of \citet{ye} can be done under a variation of Definition~\ref{ass1} where the constants $\tau_t$ and $\tau_r$ depend on the discount factor\footnote{
Define the following $\gamma$-discounted variation of $\tau^\pi(i,t)$:
$
 \tau_\gamma^\pi(i,t) = \E \left[ \sum_{k=0}^{t-1} \gamma^k {\textbf{1}}_{i_t=i}~|~i_0 \sim U,~ a_t=\pi(i_t)\right] = \sum_{k=0}^{t-1} \gamma^k \P(i_t=i~|~i_0 \sim U,~ a_t=\pi(i_t))
$
 and $\tau_\gamma^\pi(i)=\lim_{t \rightarrow \infty}\tau_\gamma^\pi(i,t)$. Assume that we have constants $\tau^\gamma_t,$ and $\tau^\gamma_r$ such that for every policy $\pi$,  $\tau_\gamma^\pi(i) \le \tau^\gamma_t$ if $i$ is a transient state for $\pi$, and $\frac{1}{(1-\gamma)\tau_\gamma^\pi(i)} \le \tau^\gamma_r$ if $i$ is recurrent for $\pi$. Then, one can derive a bound similar to that of Theorem~\ref{stocspi} where $\tau_t$ and $\tau_r$ are respectively replaced by $\tau_t^\gamma$ and $\frac{\tau_r^\gamma}{n}$. 
At a more technical level, our analysis begins by removing the dependency with respect to $\gamma$: Lemma~\ref{propass1}, page~\pageref{propass1}, shows that for every policy $\pi$, $\tau_\gamma^\pi(i) \le \tau_t$ if $i$ is a transient state for $\pi$, and $\frac{1}{(1-\gamma)\tau_\gamma^\pi(i)} \le n \tau_r$ if $i$ is recurrent for $\pi$ (this is where we pay the $O(n)$ term because the upper bound is $n \tau_r$ instead of $\tau^\gamma_r$); we then follow the line of arguments originally given by \citet{ye2}, though our more general setting induces a few technicalities (in particular in the second part of the proof of Lemma~\ref{mostdifficult} page~\pageref{mostdifficult}).
} $\gamma$.
\end{remark}
An immediate consequence of the above result is that \spi{ }is \emph{strongly polynomial} for sets of MDPs that are much larger than the deterministic MDPs mentioned in Theorem~\ref{detspiye}. 
\begin{corollary}
For any family of MDPs indexed by $n$ and $m$ such that $\tau_t$ and $\tau_r$ are polynomial functions of $n$ and $m$, \spi{ }terminates after a number of steps that is polynomial in $n$ and $m$.
\end{corollary}

\section{Similar results for \hpi ?}

One may then wonder whether similar results can be derived for \hpi. Unfortunately, and as briefly mentioned by \citet{ye2}, the line of analysis developed for \spi{ }does not seem to adapt easily to \hpi, because simultaneously switching several actions can interfere in a way such that the policy improvement turns out to be small.
We can be more precise on what actually breaks in the approach we have described so far. On the one hand, it is possible to write counterparts of Lemmas~\ref{spidetpart1} and~\ref{spistocpart1} for \hpi{ }(see \appdx{proofcounterparts} for proofs).
\begin{lemma}
\label{hpidetpart1}
If the MDP is deterministic, after at most $n$ iterations, either \hpi{ }finishes or a new cycle appears.
\end{lemma}
\begin{lemma}
\label{hpistocpart1}
After at most $(m-n) \lceil n^2 \tau_t \log (n^2 \tau_t)\rceil + n \lceil n^2 \tau_t \log (n^2 ) \rceil$ iterations, either \hpi{ }finishes or a new recurrent class appears.
\end{lemma}
On the other hand, we did not manage to adapt Lemma~\ref{spidetpart2} nor Lemma~\ref{spistocpart2}.
In fact, it is unlikely that a result similar to that of Lemma~\ref{spidetpart2} will be shown to hold for \hpi. In a recent deterministic example due to \citet{hansen2} to show that \hpi{ }may require at least $\Omega(n^2)$ iterations, new cycles are created every single iteration but the sequence of values satisfies\footnote{This MDP has an even number of states $n=2p$. The goal is to minimize the long term expected cost. The optimal value function satisfies $v_*(i)=-p^N$ for all $i$, with $N=p^2+p$. The policies generated by \hpi{ }have values $v_{\pi_k}(i) \in  (p^{N-k-1} , p^{N-k} )$. We deduce that for all iterations $k$ and states $i$, $\frac{v_*(i)-v_{\pi_{k+1}}(i)}{v_*(i)-v_{\pi_{k}}(i)} \ge \frac{1+p^{-k-2}}{1+p^{-k}} = 1-\frac{p^{-k}-p^{-k-2}}{1+p^{-k}} \ge 1-p^{-k}(1-p^{-2}) \ge 1-p^{-k}$. } for all iterations $k < \frac{n^2}{4}+\frac{n}{4}$ and states $i$,
\begin{align}
v_*(i)-v_{\pi_{k+1}}(i) \ge \left[1-\left( \frac{2}{n} \right)^k \right] (v_*(i)-v_{\pi_{k}}(i)).
\end{align}
Contrary to Lemma~\ref{spidetpart2}, as $k$ grows, the amount of contraction gets (exponentially) smaller and smaller. With respect to \spi, this suggests that \hpi{ }may suffer from subtle specific pathologies. In fact, the problem of determining the number of iterations required by \hpi{ } has been challenging for almost 30 years. It was originally identified as an open problem by \citet{schmitz}. In the simplest---deterministic---case, the complexity is still an open problem: the currently best-known lower bound is $O(n^2)$ \citep{hansen2}, while the best known upper bound is $O(\frac{m^n}{n})$~\cite{mansour,hollandersb}.

On the positive side, an adaptation of the line of proof we have considered so far can be carried out under the following assumption.
\begin{assumption}
\label{ass2}
The state space $X$ can be partitioned in two sets $\T$ and $\R$ such that for all policies $\pi$, the states of $\T$ are transient and those of $\R$ are recurrent.
\end{assumption}
Under this additional assumption, we can deduce the following original bounds.
\begin{theorem}[Proof in \appdx{proofstocpi}]
\label{stocpi}
If the MDP satisfies Assumption~\ref{ass2}, then \hpi{ }and \spi{ }terminate after at most
$$
(m-n)\left( \lceil n\tau_r \log n^2 \tau_r \rceil + \lceil n^2 \tau_t \log n^2 \tau_t \rceil \right) ~=~ \tilde O(mn (n^2\tau_t + n \tau_r))
$$
iterations.
\end{theorem}
It should however be noted that Assumption~\ref{ass2} is rather
restrictive. It implies that the algorithms converge on the recurrent
states independently of the transient states, and thus the analysis
can be decomposed in two phases: 1) the convergence on recurrent
states and then 2) the convergence on transient states (given that
recurrent states do not change anymore). The analysis of the first
phase (convergence on recurrent states) is greatly facilitated by the
fact that in this case, a new recurrent class appears every single
iteration (this is in contrast with Lemmas~\ref{spidetpart1},
\ref{spistocpart1}, \ref{hpidetpart1} and \ref{hpistocpart1} that were
designed to show under which conditions cycles and recurrent classes
are created). Furthermore, the analysis of the second phase
(convergence on transient states) is similar to that of the discounted
case of Theorems~\ref{dischpi} and \ref{discspi2}.  In other words, this last result sheds some light on the practical efficiency of \hpi{
}and \spi, and a general analysis of \hpi{ }is still largely open, and
constitutes intriguing future work.

~\\

The following sections contains detailed proofs of Lemmas~\ref{hpicontraction} and~\ref{spicontraction}, Theorems~\ref{dischpi}, \ref{discspi2}, and~\ref{stocspi}, Lemmas~\ref{hpidetpart1} and~\ref{hpistocpart1}, and finally Theorem \ref{stocpi}.
Before we start, we  provide a particularly  useful identity relating the difference between the values of two policies $\pi$ and $\pi'$  and the relative advantage $a_\pi^{\pi'}$.
\begin{lemma}
\label{id}
For all pairs of policies $\pi$ and $\pi'$,
\begin{align}
v_{\pi'}-v_{\pi} = (I-\gamma P_{\pi'})^{-1} a_\pi^{\pi'} = (I-\gamma P_{\pi})^{-1}(-a_{\pi'}^{\pi}).
\end{align}
\end{lemma}
\beginproof
This first identity follows from simple linear algebra arguments:
\begin{align}
v_{\pi'}-v_{\pi} & = (I-\gamma P_{\pi'})^{-1} r_{\pi'} - v_\pi  & \{ v_{\pi'}=T_{\pi'} v_{\pi'} ~\Leftrightarrow~ v_{\pi'}=(I-\gamma P_{\pi'})^{-1}r_{\pi'}\} \\
& = (I-\gamma P_{\pi'})^{-1} (r_{\pi'} + \gamma P_{\pi'} v_\pi-v_\pi) \\
& = (I-\gamma P_{\pi'})^{-1} (T_{\pi'}v_\pi-v_\pi).
\end{align}
The second identity follows by symmetry.
\endproofb
We will repeatedly use the following property: since for any policy $\pi$, the matrix $(1-\gamma)(I-\gamma P)^{-1}=(1-\gamma)\sum_{t=0}^\infty (\gamma P_\pi)^t$ is a stochastic matrix (as a mixture of stochastic matrices), then 
$$
\|(I-\gamma P)^{-1}\|_\infty=\frac{1}{1-\gamma},
$$
where $\|\cdot\|_\infty$ is the natural induced max-norm on matrices.
Finally, for any vector/matrix $A$ and any number $\lambda$, we shall use the notation ``$A \ge \lambda$'' (respectively ``$A \le \lambda$'') for denoting the fact that ``all the coefficients of $A$ are greater or equal to (respectively smaller or equal to) $\lambda$''.

\section{Contraction property for \hpi{ }(Proof of Lemma~\ref{hpicontraction})}

\label{proofhpicontraction}

For any $k$, we have
\begin{align}
v_{\pi_*}-v_{\pi_k} & = T_{\pi_*}v_{\pi_*} - T_{\pi_*}v_{\pi_{k-1}} + T_{\pi_*}v_{\pi_{k-1}} - T_{\pi_k}v_{\pi_{k-1}} + T_{\pi_k}v_{\pi_{k-1}} - T_{\pi_k}v_{\pi_k}  & \{ \forall \pi,~ T_\pi v_\pi=v_\pi \}  \\
& \le \gamma P_{\pi_*}(v_{\pi_*} - v_{\pi_{k-1}}) + \gamma P_{\pi_k}( v_{\pi_{k-1}}-v_{\pi_k}) & \{ T_{\pi_*}v_{\pi_{k-1}} \le T_{\pi_k}v_{\pi_{k-1}} \} \\
& \le \gamma P_{\pi_*}(v_{\pi_*}-v_{\pi_{k-1}}). & \{ \mbox{Lemma~\ref{piimp} and }P_{\pi_k} \ge 0 \}
\end{align}
Since $v_{\pi_*}-v_{\pi_k}$ is non-negative, we can take the max-norm and get:
\begin{align}
\|v_{\pi_*}-v_{\pi_k}\|_\infty \le \gamma \|v_{\pi_*}-v_{\pi_{k-1}}\|_\infty.
\end{align}

\section{Contraction property for \spi{ }(Proof of Lemma~\ref{spicontraction}) }

\label{proofspicontraction}

The proof we provide here is very close to the one given by \citet{ye}. We provide it here for completeness, and also because
it resembles the proofs we will provide for the bounds that are independent of $\gamma$.

On the one hand, using Lemma~\ref{id}, we have for any $k$:
\begin{align}
v_{\pi_{k+1}} - v_{\pi_k} & = (I-\gamma P_{\pi_{k+1}})^{-1} a_{\pi_k}^{\pi_{k+1}} \\
& \ge a_{\pi_k}^{\pi_{k+1}}, & \{ (I-\gamma P_{\pi_{k+1}})^{-1} - I \ge 0 \mbox{ and } a_{\pi_k}^{\pi_{k+1}} \ge 0 \}
\end{align}
which implies, by left multiplying by the vector $\1^T$,  that
\begin{align}
\1^T ( v_{\pi_{k+1}} - v_{\pi_k}) \ge \1^T a_{\pi_k}^{\pi_{k+1}}. \label{spieq1}
\end{align}
On the other hand, we have:
\begin{align}
 v_{\pi_*}-v_{\pi_k} & =  (I-\gamma P_{\pi_*})^{-1} a_{\pi_k}^{\pi_*}  & \{ \mbox{Lemma~\ref{id}} \}\\
& \le \frac{1}{1-\gamma} \max_s a_{\pi_k}^{\pi_{k+1}}(s) & \{ \|(I-\gamma P_{\pi_*})^{-1}  \|_\infty=\frac{1}{1-\gamma} \mbox{ and } \max_s a_{\pi_k}^{\pi_{k+1}}(s) = \max_{s,\pi} a_{\pi_k}^{\pi}(s)\ge 0\} \\
& \le  \frac{1}{1-\gamma} \1^T a_{\pi_k}^{\pi_{k+1}}, & \{ \forall x\ge 0,~\max_s x(s) \le \1^T x\}
\end{align}
which implies that
\begin{align}
\1^T  a_{\pi_k}^{\pi_{k+1}} & \ge (1-\gamma)\| v_{\pi_*}-v_{\pi_k} \|_\infty\\
& \ge  \frac{1-\gamma}{n} \1^T ( v_{\pi_*}-v_{\pi_k} ).  & \{\forall x,~\1^T x \le n\|x\|_\infty \}\label{spieq2}
\end{align}
Combining Equations~\eqref{spieq1} and~\eqref{spieq2}, we get:
\begin{align}
\1^T (v_{\pi_*}-v_{\pi_{k+1}}) 
 = \1^T(v_{\pi_*}-v_{\pi_k}) - \1^T(v_{\pi_{k+1}} - v_{\pi_k}) 
& \le \1^T(v_{\pi_*}-v_{\pi_k}) - \frac{1-\gamma}{n} \1^T (v_{\pi_*}-v_{\pi_k})\\
& = \left(1- \frac{1-\gamma}{n}\right) \1^T (v_{\pi_*}-v_{\pi_k}).
\end{align}

\section{A bound for \hpi{ }when $\gamma<1$ (Proof of Theorem~\ref{dischpi})}

\label{proofdischpi}

Although the overall line or arguments follows from those given originally by \citet{ye} and
adapted by \citet{hansen}, our proof is slightly more direct and leads to a better result.

For any $k$, we have:
\begin{align}
-a_{\pi_*}^{\pi_k}
& =  (I-\gamma P_{\pi_k})(v_*-v_{\pi_k}) & \{ \mbox{Lemma~\ref{id}} \}\\
& \le v_*-v_{\pi_k}. & \{ v_*-v_{\pi_k} \ge 0\mbox{ and }P_{\pi_k}\ge 0 \}
\end{align}
By the optimality of $\pi_*$,  $-a_{\pi_*}^{\pi_k}$ is non-negative, and we can take the max-norm:
\begin{align}
\|a_{\pi_*}^{\pi_k}\|_\infty &\le \| v_*-v_{\pi_k} \|_\infty \\
& \le \gamma^k \| v_{\pi_*}-v_{\pi_{0}}\|_\infty & \mbox{\{Lemma~\ref{hpicontraction}\}} \\
& =  \gamma^k \| (I-\gamma P_{\pi_0})^{-1} (-a_{\pi_*}^{\pi_0}) \|_\infty& \mbox{\{Lemma~\ref{id}\}} \\
& \le \frac{\gamma^k}{1-\gamma} \| a_{\pi_*}^{\pi_0} \|_\infty. & \{ \|(I-\gamma P_{\pi_0})^{-1}\|_\infty=\frac{1}{1-\gamma} \} 
\end{align}
By definition of the max-norm, and as $a_{\pi_*}^{\pi_0} \le 0$ (using again the fact that $\pi_*$ is optimal), there exists a state $s_0$ such that $-a_{\pi_*}^{\pi_0} (s_0)=\|a_{\pi_*}^{\pi_0}\|_\infty$. We deduce that for all $k$, 
\begin{align}
-a_{\pi_*}^{\pi_k}(s_0) &\le \| a_{\pi_*}^{\pi_k}\|_\infty 
 \le \frac{\gamma^k}{1-\gamma} \|a_{\pi_*}^{\pi_0}\|_\infty 
 =  \frac{\gamma^k}{1-\gamma} (-a_{\pi_*}^{\pi_0} (s_0)).
\end{align}
As a consequence, the action $\pi_k(s_0)$ must be different from $\pi_0(s_0)$ when $\frac{\gamma^k}{1-\gamma} < 1$, that is for all values of $k$ satisfying 
$$
k \ge k^* 
= \left \lceil \frac{\log  \frac{1}{1-\gamma} }{1-\gamma}   \right \rceil 
> \left \lceil \frac{\log  \frac{1}{1-\gamma} }{\log\frac{1}{\gamma}}   \right \rceil .
$$
In other words, if some policy $\pi$ is not optimal, then one of its non-optimal actions will be eliminated \emph{for good} after at most $k^*$ iterations.
By repeating this argument, one can eliminate all non-optimal actions (there are at most $n-m$ of them), and the result follows.

\section{A bound for \spi{ }when $\gamma<1$ (Proof of Theorem~\ref{discspi2})}
\label{proofdiscspi2}


At each iteration $k$, let $s_k$ be the state in which an action is switched. We have (by definition of \spi):
\begin{align}
a_{\pi_k}^{\pi_{k+1}}(s_k) & = \max_{\pi,s} a_{\pi_k}^\pi (s).
\end{align}
Starting with arguments similar to those for the contraction property of \spi, we have on the one hand:
\begin{align}
v_{\pi_{k+1}}-v_{\pi_k} & = (I-\gamma P_{\pi_{k+1}})^{-1} a_{\pi_k}^{\pi_{k+1}} & \{\mbox{Lemma \ref{id}} \}  \\
& \ge a_{\pi_k}^{\pi_{k+1}},  & \{ (I-\gamma P_{\pi_{k+1}})^{-1}-I \ge 0 \mbox{ and }a_{\pi_k}^{\pi_{k+1}} \ge 0 \}
\end{align}
which implies that
\begin{align} 
v_{\pi_{k+1}}(s_k)-v_{\pi_k}(s_k) \ge a_{\pi_k}^{\pi_{k+1}}(s_k). \label{eq2}
\end{align}
On the other hand, we have:
\begin{align}
v_{\pi_*} - v_{\pi_k} & = (I-\gamma P_{\pi_*})^{-1}a^{\pi_*}_{\pi_k} & \{\mbox{Lemma \ref{id}} \} \\
& \le \frac{1}{1-\gamma}  a^{\pi_{k+1}}_{\pi_k}(s_k) &  \{ \|(I-\gamma P_{\pi_*})^{-1}  \|_\infty=\frac{1}{1-\gamma} \mbox{ and } a^{\pi_{k+1}}_{\pi_k}(s_k) = \max_{s,\pi} a_{\pi_k}^{\pi}(s) \ge 0\}
\end{align}
which implies that
\begin{align}
\|v_{\pi_*}-v_{\pi_k} \|_\infty \le \frac{1}{1-\gamma} a_{\pi_k}^{\pi_{k+1}}(s_k). \label{eq1}
\end{align}
Write $\Delta_k=v_{\pi_*}-v_{\pi_k}$. From Equations~\eqref{eq2} and \eqref{eq1}, we deduce that:
\begin{align}
\Delta_{k+1}(s_k)&  \le \Delta_k(s_k) - (1-\gamma) \|\Delta_k\|_\infty 
 = \left(1 - (1-\gamma)\frac{\|\Delta_k\|_\infty}{\Delta_k(s_k)} \right)\Delta_k(s_k).
\end{align}
This implies---since $\Delta_k(s_k) \le \|\Delta_k\|_\infty$---that
\begin{align}
\Delta_{k+1}(s_k) \le \gamma \Delta_k(s_k),
\end{align}
but also---since $\Delta_k(s_k)$ and $\Delta_{k+1}(s_k)$ are non-negative and thus $\left(1 - (1-\gamma)\frac{\|\Delta_k\|_\infty}{\Delta_k(s_k)} \right)\ge 0$---that
\begin{align}
\|\Delta_k \|_\infty \le \frac{1}{1-\gamma}\Delta_{k}(s_k). 
\end{align}
Now, write $n_k$ for the vector on the state space such that $n_k(s)$ is the number of times state $s$ has been switched until iteration $k$ (including $k$). Since by Lemma~\ref{piimp} the sequence $(\Delta_k)_{k \ge 0}$ is non-increasing, we have
\begin{align}
\|\Delta_{k}\|_\infty & \le \frac{1}{1-\gamma}\Delta_{k}(s_k)
 \le \frac{\gamma^{n_{k-1}(s_k)}}{1-\gamma} \Delta_{0}(s_k) 
 \le \frac{\gamma^{n_{k-1}(s_k)}}{1-\gamma} \| \Delta_0 \|_\infty. \label{eq3}
\end{align}
At any iteration $k$, let $s^*_k=\arg\max_s n_{k-1}(s)$ be the state in which actions have been switched the most.
Since at each iteration $k$, one of the $n$ components of $n_k$ is increased by 1, we necessarily have 
\begin{align}
n_{k-1}(s^*_k) \ge \left\lfloor \frac{k-1}{n} \right\rfloor \ge \frac{k-n}{n}. \label{eq3b}
\end{align}
Write $k^* \le k-1$ for the last iteration when the state $s^*_k$ was updated, such that we have
\begin{align}
n_{k-1}(s^*_k)=n_{k^*-1}(s_{k^*}). \label{eq4}
\end{align} 
Since $(\|\Delta_k\|_\infty)_{k \ge 0}$ is nonincreasing (using again Lemma~\ref{piimp}), we have
\begin{align}
\|\Delta_{k}\|_\infty & \le \|\Delta_{k^*}\|_\infty & \{ k^* \le k-1 \} \\
& \le \frac{\gamma^{n_{k^*-1}(s_{k^*})}}{1-\gamma} \| \Delta_0 \|_\infty & \{ \mbox{Equation~\eqref{eq3}}  \} \\
& = \frac{\gamma^{n_{k-1}(s^*_{k})}}{1-\gamma} \| \Delta_0 \|_\infty  & \{ \mbox{Equation~\eqref{eq4}} \} \\
& \le \frac{\gamma^{\frac{k-n}{n}}}{1-\gamma} \| \Delta_0 \|_\infty.  & \{ \mbox{Equation~\eqref{eq3b}~and~$x \mapsto \gamma^x$ is decreasing}  \}
\end{align}
We are now ready to finish the proof. By using arguments similar to those for \hpi, we have:
\begin{align}
\|a_{\pi_*}^{\pi_{k}}\|_\infty &\le \|\Delta_k\|_\infty 
 \le \frac{\gamma^{\frac{k-n}{n}}}{1-\gamma} \| \Delta_0 \|_\infty 
 \le \frac{\gamma^{\frac{k-n}{n}}}{(1-\gamma)^2} \|a_{\pi_*}^{\pi_{0}}\|_\infty.  
\end{align}
In particular, we can deduce from the above relation that as soon as $\frac{\gamma^{\frac{k-n}{n}}}{(1-\gamma)^2}<1$, that is for instance when
$
k>k^*=n \left( 1 + \frac{2}{1-\gamma}\log \frac{1}{1-\gamma} \right),
$ 
one of the non-optimal actions of $\pi_0$ cannot appear in $\pi_k$. Thus, every $k^*$
iterations, a non-optimal action is eliminated for good, and the result follows from the fact that there are at most $n-m$ non-optimal actions.

\section{A general bound for \spi{ }(Proof of Theorem~\ref{stocspi})}
\label{proofstocspi}

The proof we give here is strongly inspired by that for the deterministic case of \citet{ye2}: the steps (a series of lemmas) are similar. There are mainly two differences. First, our arguments are \emph{more direct} in the sense that we do not refer to linear programming, but only provide simple linear algebra arguments. Second, it is \emph{more general}: for any policy $\pi$, we consider the set of transient states (respectively recurrent classes) instead of the set of path states (respectively cycles); it slightly complicates the arguments, the most complicated extension being the second part of the proof of the forthcoming Lemma~\ref{mostdifficult}. 

Consider the vector $x_\pi=(I-\gamma P_\pi^T)^{-1}\1$ that provides a discounted measure of state visitations along the trajectories induced by a policy $\pi$ starting from the uniform distribution $U$ on the state space $X$:
\begin{align}
\forall i \in X,~~ x_\pi(i)=n \sum_{t=0}^\infty \gamma^t \P(i_t=i~|~i_0 \sim U,~ a_t=\pi(i_t)).
\end{align}
This vector plays a crucial role in the analysis. For any policy $\pi$ and state $i$, we trivially have $x_\pi(i) \in \left(1,\frac{n}{1-\gamma}\right)$.  In the case of deterministic MDPs,  \citet{ye2}) exploits the fact that $x_\pi(i)$ belongs to the set $(1,n)$ when $i$ is on path of $\pi$, while $x_\pi(i)$ belongs to the set $(\frac{1}{1-\gamma},\frac{n}{1-\gamma})$ when $i$ is on a cycle of $\pi$. Our extension of their result to the case of general (stochastic) MDPs will rely on the following result. For any policy $\pi$, we shall write $\R(\pi)$ for the set of states that are recurrent for $\pi$. 
\begin{lemma}
\label{propass1}
With the constants $\tau_t$ and $\tau_r$ of Definition~\ref{ass1}, we have for every discount factor~$\gamma$,
\begin{align}
\forall i \not\in \R(\pi),~~~1 ~\le~  & x_\pi(i) ~\le~ n \tau_t \label{prop2}  \\
\forall i \in \R(\pi),~~~ \frac{1}{\tau_r}~\le~ & (1-\gamma)x_\pi(i)   ~\le~ n. \label{prop1} 
\end{align}
\end{lemma}
\beginproof
The fact that $x_\pi(i)$ belongs to $\left(1,\frac{n}{1-\gamma}\right)$ is obvious from the definition of $x_\pi$. The upper bound on $x_\pi$ on the transient states $i$ follows from the fact that for any policy $\pi$,
\begin{align}
\tau_t(i) & \ge \lim_{t \rightarrow \infty} \tau^{\pi}(i,t) \\
& = \sum_{k=0}^{\infty} \P(i_t=i~|~i_0 \sim U,~ a_t=\pi(i_t)) \\
& \ge \sum_{k=0}^{\infty} \gamma^k \P(i_t=i~|~i_0 \sim U,~ a_t=\pi(i_t)) \\
& = \frac{1}{n} x_\pi(i).
\end{align}
Let us now consider the lower bound on $(1-\gamma)x_\pi(i)$ when $i$ is a recurrent state of some policy $\pi$. In general, the asymptotic frequency $\mu^\pi$ of $\pi$ does not necessarily satisfy ${\mu^\pi}^T P_\pi = P_\pi$ because $P_\pi$ may correspond to an aperiodic or reducible chain. To deal with this issue, we consider the
Ces\`aro mean 
\begin{align}
Q_\pi= \lim_{t \rightarrow \infty} \frac{1}{t} \sum_{k=0}^{t-1} (P_\pi)^k
\end{align}
that is well-defined~\cite[Section~3.2]{stroock}.
It can be shown~\cite[Proposition~3.5(a)]{fritz} that $Q_\pi=Q_\pi P_\pi = P_\pi Q_\pi = Q_\pi Q_\pi$. This implies in particular that
\begin{align}
(1-\gamma) Q_\pi (I-\gamma P_\pi)^{-1} = (1-\gamma)\sum_{k=0}^\infty \gamma^k Q_\pi (P_\pi)^k = (1-\gamma)\sum_{k=0}^\infty \gamma^k Q_\pi  = Q_\pi. \label{qpq}
\end{align}
Then, by using twice the fact that $\mu^\pi=\frac{1}{n}{Q_\pi}^T \1 $, we can see that for all recurrent states $i$,
\begin{align}
\frac{1}{\tau_r} & \le \mu^\pi(i) \\
& = \left[ \frac{1}{n} {Q_\pi}^T \1 \right](i) \\
& =   \left[ \frac{1}{n} (1-\gamma)(I-\gamma {P_\pi}^T)^{-1}{Q_\pi}^T \1  \right](i) & \{ \mbox{Equation~\eqref{qpq}} \} \\
& =  \left[ (1-\gamma)(I-\gamma {P_\pi}^T)^{-1} \mu^\pi \right](i) \\
& \le  \left[ (1-\gamma)(I-\gamma {P_\pi}^T)^{-1} \1 \right](i) & \{ \mu^\pi \le \1 \} \\
& = (1-\gamma)x_\pi(i).
\end{align}
\endproofb

Finally, a rewriting of Lemma~\ref{id} in terms of the vector $x_\pi$ will be useful in the following proofs: for any pair of policies $\pi$ and $\pi'$,
\begin{align}
\1^T (v_{\pi'}-v_\pi) = {x_{\pi'}}^T a^{\pi'}_\pi = {x_{\pi}}^T (-a^{\pi}_{\pi'}). \label{id2}
\end{align}

We are now ready to delve into the details of the arguments.
As mentioned before, the proof is structured in two steps: first, we will show that recurrent classes are created often; then we will show that significant progress is made every time a new recurrent class appears.

\subsection{Part 1: Recurrent classes are created often}

\begin{lemma}
\label{contractspi}
Suppose one moves from policy $\pi$ to policy $\pi'$ \emph{without creating any recurrent class}. Let $\pi_\dag$ be the final policy before either a new recurrent class appears or \spi{ }terminates. Then
\begin{align}
\1^T (v_{\pi_\dag}-v_{\pi'}) \le \left( 1 - \frac {1}{n^2 \tau_t}  \right) \1^T (v_{\pi_\dag}-v_{\pi}).
\end{align} 
\end{lemma}
\beginproof
The arguments are similar to those for the proof of Theorem~\ref{discspi2}. On the one hand, we have:
\begin{align}
\1^T (v_{\pi'}-v_\pi) \ge \1^T a^{\pi'}_{\pi}. \label{i1}
\end{align}
On the other hand, we have
\begin{align}
\1^T (v_{\pi_\dag}-v_{\pi}) & = x_{\pi_\dag}^T a^{\pi_\dag}_\pi & \{ \mbox{Equation~\eqref{id2}} \} \\
& = \sum_{s \not\in \R(\pi_\dag)} x_{\pi_\dag}(s)a^{\pi_\dag}_\pi(s) + \sum_{s \in \R(\pi_\dag)} x_{\pi_\dag}(s)a^{\pi_\dag}_\pi(s) \\
& \le n^2 \tau_t \max_{s \not\in \R(\pi_\dag)} a^{\pi_\dag}_\pi(s) + \frac{n^2}{1-\gamma} \max_{s \in \R(\pi_\dag)} a^{\pi_\dag}_\pi(s). & \{ \mbox{Equations~\eqref{prop2}-\eqref{prop1}} \}
\end{align}
Since by assumption recurrent classes of $\pi_\dag$ are also recurrent classes of $\pi$, we deduce that for all $s \in \R(\pi_\dag)$, $\pi_\dag(s)=\pi(s)$, so that  $\max_{s \in \R(\pi_\dag)} a^{\pi_\dag}_\pi(s)=0$. Thus, the second term of the above r.h.s. is null and
\begin{align}
\1^T (v_{\pi_\dag}-v_{\pi}) & \le n^2 \tau_t \max_s a^{\pi_\dag}_\pi(s)\\ 
& \le n^2 \tau_t \max_s a^{\pi'}_\pi(s) & \{\max_s T_{\pi'}v_\pi(s)=\max_{s,\tilde \pi} T_{\tilde\pi} v_\pi(s) \}\\
& \le n^2 \tau_t \1^T a^{\pi'}_\pi. \label{i2} & \{ \forall x \ge 0,~\max_s x(s) \le \1^T x\}
\end{align}
Combining Equations~\eqref{i1} and \eqref{i2}, we get:
\begin{align}
\1^T (v_{\pi_\dag}-v_{\pi'}) & = \1^T (v_{\pi_\dag}-v_{\pi}) - \1^T (v_{\pi'}-v_{\pi}) \\
& \le \left(1-\frac{1}{n^2 \tau_t} \right)  \1^T (v_{\pi_\dag}-v_{\pi}). 
\end{align}
\endproofb

\begin{lemma}
\label{mostdifficult}
While \spi{ }does not create any recurrent class nor finishes,
\begin{itemize}
\item either an action is eliminated from policies after at most $\lceil n^2 \tau_t \log(n^2 \tau_t) \rceil$ iterations,
\item or a recurrent class is broken after at most $\lceil n^2 \tau_t \log(n^2) \rceil$ iterations.
\end{itemize}
\end{lemma}
\beginproof
Let $\pi$ be the policy in some iteration. Let $\pi_\dag$ be the last policy before a new recurrent class appears, and $\pi'$ a policy generated after $k$ iterations from $\pi$. We shall prove that one of the two events stated of the lemma must happen.\\
Since 
\begin{align}
0 &\le \1^T (v_{\pi_\dag}-v_\pi) & \{ v_{\pi_\dag} \ge v_\pi \} \\
&= {x_\pi}^T (-a^\pi_{\pi_\dag}) & \{ \mbox{Equation \eqref{id2}} \} \\
&= \sum_{s \not\in \R(\pi)} x_\pi(s)  (-a^\pi_{\pi_\dag}(s)) + \sum_{C \in \R(\pi)}\sum_{s \in C} x_\pi(s)  (-a^\pi_{\pi_\dag}(s))
\end{align}
there must exist either a state $s_0 \not\in \R(\pi)$ such that
\begin{align}
x_\pi(s_0) (-a^\pi_{\pi_\dag}(s_0)) \ge \frac{1}{n} {x_\pi}^T (-a^\pi_{\pi_\dag})  \label{case1} \ge 0.
\end{align}
or a recurrent class $R_0$ such that
\begin{align}
\sum_{s \in R_0} x_\pi(s) (-a^\pi_{\pi_\dag}(s)) \ge \frac{1}{n} {x_\pi}^T (-a^\pi_{\pi_\dag})   \ge 0 \label{case2}.
\end{align}
We consider these two cases separately below.
\begin{itemize}
\item {\bf case 1:} Equation~\eqref{case1} holds for some $s_0 \not\in \R(\pi)$.
Let us prove by contradiction that for $k$ sufficiently big,  $\pi'(s_0)\neq \pi(s_0)$: let us assume that $\pi'(s_0)=\pi(s_0)$.  Then
\begin{align}
\1^T (v_{\pi_\dag}-v_{\pi'}) & \ge v_{\pi_\dag}(s_0)-v_{\pi'}(s_0) & \{v_{\pi_\dag} \ge v_{\pi'} \} \\
& = v_{\pi_\dag}(s_0)-T_{\pi'}v_{\pi'}(s_0) & \{ v_{\pi'}=T_{\pi'}v_{\pi'}\}\\
& \ge v_{\pi_\dag}(s_0)-T_{\pi'}v_{\pi_\dag}(s_0) & \{v_{\pi_\dag} \ge v_{\pi'}\} \\
& = -a^{\pi'}_{\pi_\dag}(s_0) \\
& = -a^{\pi}_{\pi_\dag}(s_0) & \{ \pi(s_0)=\pi'(s_0) \} \\
& \ge \frac{1}{n \tau_t} x_\pi(s_0)(-a^{\pi}_{\pi_\dag}(s_0)) & \{ \mbox{Equation \eqref{prop2}}\}\\
& \ge \frac{1}{n^2 \tau_t} {x_\pi}^T (-a^{\pi}_{\pi_\dag}) & \{ \mbox{Equation } \eqref{case1} \}\\
& = \frac{1}{n^2 \tau_t} \1^T (v_{\pi_\dag}-v_\pi). & \{\mbox{Equation } \eqref{id2} \}
\end{align}
If there is no recurrent class creation, the contraction property given in Lemma~\ref{contractspi} implies that if $\pi'$ is obtained after $k=\lceil n^2 \tau_t \log(n^2 \tau_t ) \rceil >  \frac{\log (n^2 \tau_t )}{\log \frac{1}{1-\frac 1 {n^2 \tau_t}}}$ iterations, then
\begin{align}
\1^T (v_{\pi_\dag}-v_{\pi'})  < \frac{1}{n^2 \tau_t} \1^T (v_{\pi_\dag}-v_\pi),
\end{align}
and we get a contradiction. As a conclusion, we necessarily have $\pi'(s_0)\neq \pi(s_0)$.
\item {\bf case 2:} Equation~\eqref{case2} holds for some $R_0$ that is a recurrent class of $\pi$. Let us prove by contradiction that for $k$ sufficiently big, $R_0$ cannot be a recurrent class of  $\pi'$: let us thus assume that $R_0$ is a recurrent class of $\pi'$.
Write ${\cal T}$ for the set of states that are transient for $\pi$ (formally, ${\cal T}=X\backslash \R(\pi)$).
For any subset $Y$ of the state space $X$, write $P_{\pi}^Y$ for the stochastic matrix of which the $i^{th}$ row is equal to that of $P_\pi$ if $i \in Y$, and is $0$ otherwise, and write $\1_Y$ the vectors of which the $i^{th}$ component is equal to 1 if $i \in Y$ and 0 otherwise.

Using the fact that $P_\pi^{R_0}P_\pi^{{\cal T}}=0$, one can first observe that
\begin{align}
(I-\gamma P_\pi^{R_0})(I-\gamma P_\pi^{{\cal T}})=I-\gamma (P_\pi^{R_0}+P_\pi^{{\cal T}}),
\end{align}
from which we can deduce that
\begin{align}
\forall s \in R_0,~\left[ {\1_{{\cal T} \cup R_0}}^T (I-\gamma P_\pi)^{-1}\right] (s) & = \left[{\1_{{\cal T} \cup R_0}}^T (I-\gamma (P_\pi^{R_0}+P_\pi^{{\cal T}}))^{-1} \right] (s) \\
& = \left[ {\1_{{\cal T} \cup R_0}}^T (I-\gamma P_\pi^{{\cal T}})^{-1}(I-\gamma P_\pi^{R_0})^{-1}\right] (s). \label{trick}
\end{align}

Also, let $s$ be an arbitrary state and $s'$ be a state of $R_0$. Since $\1_s^T (P_\pi^{\cal T})^k (s')$ is the probability that the chain starting in $s$ reaches $s'$ for the first time after $k$ iterations, then
$$
\1_s^T (I-\gamma {P_\pi^{{\cal T}}})^{-1} (s') \le \sum_{t=0}^\infty  \1_s^T(P_\pi^{\cal T})^{i} (s') \le 1.
$$
and therefore,
\begin{align}
\forall s' \in R_0,~~~{\1_{{\cal T} \cup R_0}}^T(I-\gamma {P_\pi^{{\cal T}}})^{-1} (s')  \le n. \label{HP2}
\end{align}
Writing $\delta$ for the vector that equals $-a^{\pi}_{\pi_\dag}$ on $R_0$ and that is null everywhere else, we have
\begin{align}
& \sum_{s \in R_0} x_\pi(s)(-a^{\pi}_{\pi_\dag}(s)) \\
& = \sum_{s \in R_0}  [(I-\gamma P_\pi^T)^{-1} \1] (s) \delta(s) \\
& = \sum_{s \in R_0}  [(I-\gamma P_\pi^T)^{-1} {\1_{{\cal T} \cup R_0}}] (s) \delta(s) & \left\{  \forall s \in R_0,~[(I-\gamma P_\pi^T)^{-1} \1_{X \backslash ({\cal T}\cup R_0)}] (s)=0  \right\}~~~ \\
& = \sum_{s} [(I-\gamma P_\pi^T)^{-1} {\1_{{\cal T} \cup R_0}}] (s) \delta(s) & \{ \forall s \not\in R_0,~\delta(s)=0\}~~~\\
& = {\1_{{\cal T} \cup R_0}}^T (I-\gamma P_\pi)^{-1} \delta \\
& = {\1_{{\cal T} \cup R_0}}^T (I-\gamma P_\pi^{{\cal T}})^{-1}(I-\gamma P_\pi^{R_0})^{-1} \delta & \{\mbox{Equation~\eqref{trick}}\}~~~ \\
& = \sum_{s} [ (I-\gamma {P_\pi^{{\cal T}}}^T)^{-1}{\1_{{\cal T} \cup R_0}}](s) [(I-\gamma P_\pi^{R_0})^{-1} \delta](s) \\
& = \sum_{s \in R_0} [ (I-\gamma {P_\pi^{{\cal T}}}^T)^{-1}{\1_{{\cal T} \cup R_0}}](s) [(I-\gamma P_\pi^{R_0})^{-1} \delta](s) & \{ \forall s \not\in R_0,~\delta(s)=0\}~~~ \\
& = \sum_{s \in R_0} [ (I-\gamma {P_\pi^{{\cal T}}}^T)^{-1}{\1_{{\cal T} \cup R_0}}](s) (v_{\pi_\dag}(s)-v_\pi(s)) & \{ \mbox{Lemma~\ref{id}} \}~~~ \\
& \le n {\1_{R_0}}^T (v_{\pi_\dag}-v_{\pi}). & \{ \mbox{Equation \eqref{HP2}} \}~~~ \label{HP3}
\end{align}
We assumed that $R_0$ is also a recurrent class of $\pi'$, which implies ${\1_{R_0}}^T v_\pi = {\1_{R_0}}^T v_{\pi'}$, and
\begin{align}
\1^T (v_{\pi_\dag}-v_{\pi'}) & \ge {\1_{R_0}}^T( v_{\pi_\dag} -v_{\pi'}) & \{ v_{\pi_\dag} \ge v_{\pi'} \} \\
& = {\1_{R_0}}^T( v_{\pi_\dag} -v_{\pi} ) & \{ {\1_{R_0}}^T v_\pi = {\1_{R_0}}^T v_{\pi'} \} \\
 & \ge \frac{1}{n}\sum_{s \in R_0} x_\pi(s) (-a^{\pi}_{\pi_\dag}(s)) & \{ \mbox{Equation \eqref{HP3}}\}\\
& \ge \frac{1}{n^2} {x_\pi}^T (-a^{\pi}_{\pi_\dag}) & \{\mbox{Equation } \eqref{case2} \}\\
& = \frac{1}{n^2} \1^T (v_{\pi_\dag}-v_\pi). & \{ \mbox{Equation \eqref{id2}} \}
\end{align}
If there is no recurrent class creation, the contraction property given in Lemma~\ref{contractspi} implies that if $\pi'$ is obtained after $k=\lceil n^2 \tau_t \log(n^2) \rceil >  \frac{\log (n^2)}{\log \frac{1}{1-\frac 1 {n^2 \tau_t}}}$ iterations, then 
\begin{align}
\1^T (v_{\pi_\dag}-v_{\pi'})  < \frac{1}{n^2} \1^T (v_{\pi_\dag}-v_\pi),
\end{align}
and thus we get a contradiction. As a conclusion, $R_0$ cannot be a recurrent class of $\pi'$.
\end{itemize}
\endproofb

A direct consequence of the above result is Lemma~\ref{spistocpart1} that we originally stated on page~\pageref{spistocpart1}, and that we restate for clarity.
\setcounter{lemma}{5} 
\begin{lemma}
After at most $(m-n) \lceil n^2 \tau_t \log (n^2 \tau_t)\rceil + n \lceil n^2 \tau_t \log (n^2 ) \rceil $ iterations, either \spi{ }finishes or a new recurrent class appears.
\end{lemma}
\beginproof
Before a recurrent class is created, at most $n$ recurrent classes need to be broken and $(m-n)$ actions to be eliminated, and the time required by these events is bounded thanks to the previous lemma.
\endproofb

\subsection{Part 2: A new recurrent class implies a significant step towards the optimal value}

We now proceed to the second part of the proof, and begin by proving Lemma~\ref{spistocpart2} (originally stated page~\pageref{spistocpart2}).
\begin{lemma}
When \spi{ }moves from $\pi$ to $\pi'$ where $\pi'$ involves a new recurrent class, we have
\begin{align}
\1^T(v_{\pi_*}-v_{\pi'})  \le \left( 1-\frac{1}{n \tau_r} \right) \1^T(v_{\pi_*}-v_{\pi}).
\end{align}
\end{lemma}
\setcounter{lemma}{13}
\beginproof
Let $s_0$ be the state such that $\pi'(s_0) \neq \pi(s_0)$. On the one hand, since $\pi'$ contains a new recurrent class $R$ (necessarily containing $s_0$), we have
\begin{align}
\1^T (v_{\pi'}-v_\pi) 
& = {x_{\pi'}}^T a^{\pi'}_\pi & \{ \mbox{Equation \eqref{id2}} \}\\
& = x_{\pi'}(s_0) a_\pi(s_0) & \{\mbox{\spi{ }switches 1 action and $a_\pi(s_0)=a_\pi^{\pi'}(s_0)$} \} \\
& \ge \frac{1}{(1-\gamma)\tau_r} a_\pi(s_0). & \{ \mbox{Equation~\eqref{prop1} with }s_0 \in \R(\pi')\} \label{i3}
\end{align}
On the other hand, 
\begin{align}
\forall s,~~~v_{\pi_*}(s)-v_{\pi}(s) & = [(I-\gamma P_{\pi_*})^{-1} a^{\pi_*}_\pi](s)  & \{ \mbox{Lemma~\ref{id}} \} \\
& \le \frac{1}{1-\gamma}  a_\pi(s_0). & \{ \|(I-\gamma P_{\pi_*})^{-1}\|_\infty \le \frac{1}{1-\gamma} \mbox{ and } a_\pi(s_0) = \max_{s,\tilde \pi}a^{\tilde\pi}_\pi(s) \ge 0 \} \label{i4}\\
\end{align}
Combining these two observations, we obtain
\begin{align}
\1^T(v_{\pi_*}-v_{\pi'}) & = \1^T(v_{\pi_*}-v_{\pi})-\1^T(v_{\pi'}-v_{\pi}) \\
& \le  \1^T(v_{\pi_*}-v_{\pi}) - \frac{1}{(1-\gamma)\tau_r} a_\pi(s_0) & \{ \mbox{Equation \eqref{i3}} \}\\
& \le  \1^T(v_{\pi_*}-v_{\pi}) - \frac{1}{\tau_r} \max_s v_{\pi*}(s)-v_{\pi'}(s)  & \{ \mbox{Equation \eqref{i4}} \}\\
& \le \left( 1-\frac 1 {n \tau_r} \right) \1^T(v_{\pi_*}-v_{\pi}). & \{ \forall x,~ \frac{1}{n}\1^T x \le \max_s x(s) \}
\end{align}
\endproofb

\begin{lemma}
\label{spistocelim}
While \spi{ }does not terminate, 
\begin{itemize}
\item either some non-optimal action is eliminated from recurrent states after at most $\lceil n \tau_r \log (n^2 \tau_r) \rceil$ recurrent class creations,
\item or some non-optimal action is eliminated from policies after at most $\lceil n \tau_r \log (n^2 \tau_t) \rceil$ recurrent class creations.
\end{itemize}
\end{lemma}
\beginproof
Let $\pi$ be the policy in some iteration and $\pi'$ the policy generated after $k$ iterations from $\pi$ (without loss of generality we assume $\pi'\neq \pi_*$).
Let $s_0=\arg\max_s x_\pi (s)(-a^\pi_{\pi_*}(s))$. We have 
\begin{align}
x_\pi (s_0)( -a^\pi_{\pi_*}(s_0)) & \ge \frac{1}{n} {x_\pi}^T (-a^\pi_{\pi_*})  & \{ \forall x,~\1^T x \le n \max_s x(s) \}\\
& = \frac{1}{n} \1^T (v_{\pi_*}-v_\pi). & \{ \mbox{Equation \eqref{id2}}\}  \label{fact2}
\end{align}
We now consider two cases, respectively corresponding to $s_0 \not\in \R(\pi)$ or $s_0 \in \R(\pi)$.
\begin{itemize}
\item {\bf case 1:} $s_0 \not\in \R(\pi)$.
Let us prove by contradiction that $\pi'(s_0)\neq \pi(s_0)$ if $k$ is sufficiently large: let us assume that $\pi'(s_0)=\pi(s_0)$. Then, by using repeatedly the fact that for all $\tilde\pi,~ a^{\tilde\pi}_{\pi_*} \le 0$ (by definition of the optimal policy~$\pi_*$), we have:
\begin{align}
\1^T (v_{\pi_*}-v_{\pi'}) & = {x_{\pi'}}^T (-a^{\pi'}_{\pi_*}) & \{ \mbox{Equation \eqref{id2}}  \}\\
& \ge x_{\pi'}(s_0) (-a^{\pi'}_{\pi_*} (s_0))  \\
&\ge -a^{\pi'}_{\pi_*} (s_0) & \{ x_{\pi'}(s_0) \ge 1  \}\\
& = -a^{\pi}_{\pi_*}(s_0) & \{ \pi(s_0)=\pi'(s_0) \}\\
& \ge \frac{1}{n \tau_t}x_\pi(s_0)(-a^{\pi}_{\pi_*}(s_0))  & \{ \mbox{Equation \eqref{prop2}} \} \\
& \ge \frac{1}{n^2 \tau_t} \1^T (v_{\pi_*}-v_\pi). & \{ \mbox{Equation }\eqref{fact2} \}
\end{align}
After $k=\lceil n \tau_r \log n^2 \tau_t \rceil > \frac{\log n^2 \tau_t }{\log \frac{1}{1-\frac 1 {n \tau_r}}}$ recurrent classes are created, we have by the contraction property of Lemma~\ref{spistocpart2} that
\begin{align}
\1^T (v_{\pi_*}-v_{\pi'})  < \frac{1}{n^2 \tau_t} \1^T (v_{\pi_*}-v_\pi)
\end{align}
and we get a contradiction. As a conclusion, we have $\pi'(s_0) \neq \pi(s_0)$. 
\item {\bf case 2:} $s_0 \in \R(\pi)$.
Let us prove by contradiction that $\pi'(s_0) \neq \pi(s_0)$ if $s_0$ is recurrent for $\pi'$ and $k$ is sufficiently large: let us assume that $\pi'(s_0)=\pi(s_0)$ and $s_0 \in \R(\pi')$. Then, by using again the fact that for all $\tilde\pi,~ a^{\tilde\pi}_{\pi_*} \le 0$, we have:
\begin{align}
\1^T (v_{\pi_*}-v_{\pi'}) &= {x_{\pi'}}^T (-a^{\pi'}_{\pi_*}) & \{ \mbox{Equation \eqref{id2}} \}\\
& = \sum_{s} x_{\pi'}(s)(-a^{\pi'}_{\pi_*}(s)) \\
& \ge \sum_{s \in R_0} x_{\pi'}(s)(-a^{\pi'}_{\pi_*}(s)) \\
& \ge \frac{1}{(1-\gamma)\tau_r} \sum_{s \in R_0}(-a^{\pi'}_{\pi_*}(s))  & \{ \mbox{Equation~\eqref{prop1}} \} \\
& \ge \frac{1}{(1-\gamma)\tau_r} (-a^{\pi'}_{\pi_*}(s_0)) \\
& = \frac{1}{(1-\gamma)\tau_r} (-a^{\pi}_{\pi_*}(s_0)) & \{ \pi(s_0)=\pi'(s_0) \}\\
& \ge \frac{1}{n \tau_r}x_\pi (s_0)(-a^{\pi}_{\pi_*}(s_0)) & \{ x_\pi(s_0) \le \frac{n}{1-\gamma} \}\\
& \ge \frac{1}{n^2 \tau_r} \1^T (v_{\pi_*}-v_\pi). & \{ \mbox{Equation }\eqref{fact2} \}
\end{align}
After $k=\lceil n \tau_r \log n^2 \tau_r \rceil > \frac{\log n^2 \tau_r }{\log \frac{1}{1-\frac 1 {n \tau_r}}}$ new recurrent classes are created, we have by the contraction property of Lemma~\ref{spistocpart2} that
\begin{align}
\1^T (v_{\pi_*}-v_{\pi'})  < \frac{1}{n^2 \tau_r} \1^T (v_{\pi_*}-v_\pi),
\end{align}
and we get a contradiction. As a conclusion, we know that $\pi'(s_0) \neq \pi(s_0)$ if $s_0$ is recurrent for $\pi'$.
\end{itemize}
\endproofb

We are ready to conclude:
At most, the $(m-n)$ non-optimal actions may need to be eliminated from all states; in addition, all actions may need to be eliminated from recurrent states (some optimal actions may only be used at transient states and thus also need to be eliminated from recurrent states). Overall, convergence can thus be obtained after at most a total of $m \lceil n \tau_r \log (n^2 \tau_r) \rceil + (m-n) \lceil n \tau_r \log (n^2 \tau_t) \rceil$ recurrent class creations. The result follows from the fact that each class creation requires at most $(m-n) \lceil n^2 \tau_t \log (n^2 \tau_t)\rceil + n \lceil n^2 \tau_t \log (n^2 ) \rceil$ iterations (cf. Lemma~\ref{spistocpart1}).

\section{Cycle and recurrent classes creations for \hpi{ }(Proofs of Lemmas~\ref{hpidetpart1} and~\ref{hpistocpart1})}

\label{proofcounterparts}

\setcounter{lemma}{7}
\begin{lemma}
If the MDP is deterministic, after at most $n$ iterations, either \hpi{ }finishes or a new cycle appears.
\end{lemma}
\beginproof
Consider a sequence of $l$ generated policies $\pi_1,\cdots,\pi_l$ from an initial policy $\pi_0$ such that no new cycle appears. By induction, we have
\begin{align}
v_{\pi_l}-v_{\pi_k} &= T_{\pi_l}v_{\pi_l} - T_{\pi_l}{v_{\pi_{k-1}}} + T_{\pi_l}{v_{\pi_{k-1}}} - T_{\pi_k}{v_{\pi_{k-1}}} + T_{\pi_k}{v_{\pi_{k-1}}} - T_{\pi_k}{v_{\pi_{k}}} & \{ \forall \pi,~ T_\pi v_\pi=v_\pi \}~~~  \\
& \le \gamma P_{\pi_l}(v_{\pi_l} - v_{\pi_{k-1}}) + \gamma P_{\pi_k}( v_{\pi_{k-1}}-v_{\pi_k}) & \{ T_{\pi_l}v_{\pi_{k-1}} \le T_{\pi_k}v_{\pi_{k-1}} \}~~~ \\
& \le \gamma P_{\pi_l}(v_{\pi_l}-v_{\pi_{k-1}}) & \{ \mbox{Lemma~\ref{piimp} and }P_{\pi_k} \ge 0 \}~~~ \\
& \le (\gamma P_{\pi_l})^{k} (v_{\pi_l}-v_{\pi_0}). & \{ \mbox{By induction on $k$}  \}~~~ \label{rec}
\end{align}
Since the MDP is deterministic and has $n$ states,  $(P_{\pi_l})^n$ will  only have non-zero values on columns that correspond to $\R(\pi_l)$. Furthermore, since no cycle is created, $\R(\pi_l) \subset \R(\pi_0)$, which implies that $v_{\pi_l}(s)-v_{\pi_0}(s)=0$ for all $s \in \R(\pi_l)$. As a consequence, we have $(P_{\pi_l})^n (v_{\pi_l}-v_{\pi_0})=0$. By Equation~\eqref{rec}, this implies that $v_{\pi_l}=v_{\pi_n}$. If $l>n$, then \hpi{ }must have terminated.
\endproofb

\begin{lemma}
After at most $(m-n) \lceil n^2 \tau_t \log (n^2 \tau_t)\rceil + n \lceil n^2 \tau_t \log (n^2 ) \rceil$ iterations, either \hpi{ }finishes or a new recurrent class appears.
\end{lemma}
\setcounter{lemma}{17}
\beginproof
A close examination of the proof of Lemma~\ref{spistocpart1}, originally designed for \spi, shows that it applies to \hpi{ }without any modification.
\endproofb

\section{A bound for \hpi{ }and \spi{ } under Assumption~\ref{ass2} (Proof of Theorem~\ref{stocpi})}

\label{proofstocpi}

We here consider that the state space is decomposed into 2 sets: $\T$
is the set of states that are transient under all policies, and $\R$
is the set of states that are recurrent under all policies. From this
assumption, it can be seen that when running \hpi{ }or \spi, the values and actions chosen on $\T$ have no
influence on the evolution of the values and policies on $\R$. So we will study the convergence of both algorithms in two steps: we will first bound the number of iterations to converge on $\R$; we will then add the number of iterations for converging on $\T$ given that convergence has occurred on $\R$.

\paragraph{Convergence on the set $\R$ of recurrent states:}
Without loss of generality, we consider here that the state space is only made of the set of recurrent states.

First consider \spi. If all states are recurrent, new recurrent classes are created at every iteration, and Lemma~\ref{spistocpart2} holds.  Then, in a way similar to the proof of Lemma~\ref{spistocelim}, it can be shown that every $\lceil n\tau_r \log n^2 \tau_r \rceil$ iterations, a non-optimal action can be eliminated. As there are at most $(m-n)$ non-optimal actions, we deduce that \spi{ }converges in at most $(m-n) \lceil n \tau_r \log n^2 \tau_r \rceil$ iterations on $\R$.

Consider now \hpi. We can prove the following lemma.
\setcounter{lemma}{9}
\begin{lemma}
If the MDP satisfies Assumption~\ref{ass1} and all states are recurrent under all policies, \hpi{ }generates policies $(\pi_k)_{k \ge 0}$ that satisfy:
\begin{align}
\1^T(v_{\pi_*}-v_{\pi_{k+1}})  \le \left( 1-\frac{1}{n \tau_r} \right) \1^T(v_{\pi_*}-v_{\pi_k}).
\end{align}
\end{lemma}
\beginproof
On the one hand, we have
\begin{align}
\1^T (v_{\pi_{k+1}}-v_{\pi_k}) &= {x_{\pi_{k+1}}}^T a^{\pi_{k+1}}_{\pi_k} & \{\mbox{Equation \eqref{id2}}\} \\
& = {x_{\pi_{k+1}}}^T a_{\pi_k} & \{ a^{\pi_{k+1}}_{\pi_k}=a_{\pi_k} \} \\
& \ge \frac{1}{(1-\gamma)\tau_r} \1^T a_{\pi_k} & \{\mbox{Equation \eqref{prop1} and all states are recurrent}\} \\
& \ge \frac{1}{(1-\gamma)\tau_r} \| a_{\pi_k} \|_\infty. & \{\forall x \ge 0, \1^T x \ge \|x\|_\infty \} \label{coco1}
\end{align}
On the other hand, 
\begin{align}
\1^T (v_{\pi_{*}}-v_{\pi_k}) & = {x_{\pi_{*}}}^T a^{\pi_*}_{\pi_k} & \{\mbox{Equation \eqref{id2}} \} \\
& \le {x_{\pi_{*}}}^T a_{\pi_k} & \{ a_{\pi_k} \ge a_{\pi_k}^{\pi_*} \} \\
& \le \frac{n}{1-\gamma} \| a_{\pi_k} \|_\infty. & \{ \sum_i x_{\pi_{*}}(i) \le \frac{n}{1-\gamma} \mbox{ and } a_{\pi_k}\ge 0\} \label{coco2}
\end{align}
By combining Equations~\eqref{coco1} and \eqref{coco2}, we obtain:
\begin{align}
\1^T (v_{\pi_{*}}-v_{\pi_{k+1}}) & = \1^T (v_{\pi_{*}}-v_{\pi_k}) - \1^T (v_{\pi_{k+1}}-v_{\pi_k}) \\
& \le \left( 1-\frac{1}{n\tau_r}\right) \1^T (v_{\pi_{*}}-v_{\pi_k}). 
\end{align}
\endproofb
Then, similarly to \spi, we can prove that after every $\lceil n \tau_r \log n^2 \tau_r \rceil$ iterations a non-optimal action must be eliminated. And as there are at most $(m-n)$ non-optimal actions, we deduce that \hpi{ }converges in at most $(m-n)\lceil n\tau_r \log n^2 \tau_r \rceil$ iterations on $\R$.

\paragraph{Convergence on the set $\T$ of transient states:}
Consider now that convergence has occurred on the recurrent states $\R$. A simple variation of the proof of Lemma~\ref{spistocpart1}/Lemma~\ref{hpistocpart1} (where we use the fact that we don't need to consider the events where recurrent classes are broken since recurrent classes do not evolve anymore) allows us to show that the extra number of iterations for both algorithms to converge on the transient states is at most $(m-n))\lceil n^2 \tau_t \log n^2 \tau_t \rceil$, and the result follows.

\section*{Acknowledgements.}
I would like to thank Ian Post for exchanges about the proof in
\cite{ye2}, Thomas Dueholm Hansen for noticing a flaw in a claimed
result for deterministic MDPs in an earlier version, Romain Aza{\"i}s
for the reference on the Cesaro mean of stochastic matrices, and the
reviewers and editor for their very careful feedback, who helped
improve the paper overall, and the proof of Lemma 13 in particular.

\bibliographystyle{natbib2}
\bibliography{biblio.bib} 

\end{document}